\documentclass[a4paper,11pt]{article}

\usepackage{amsfonts}
\usepackage{epsf,latexsym,amsfonts,amsbsy,mathrsfs}
\usepackage{cite}
\usepackage{amsmath, amssymb, amsthm,amscd,amsxtra}
\usepackage[]{graphicx,subfigure}
\usepackage{float}
\usepackage{multirow}
\usepackage{bbm}
\usepackage{stmaryrd}
\usepackage{epsfig}

\newtheorem{theorem}{Theorem}[section]
\newtheorem{lemma}[theorem]{Lemma}

\theoremstyle{definition}
\newtheorem{definition}[theorem]{Definition}
\newtheorem{example}[theorem]{Example}

\theoremstyle{remark}
\newtheorem{remark}[theorem]{Remark}

\numberwithin{equation}{section}

\newcommand{\eproof}{$\quad \Box$}

\title{Domain Decomposition Methods for Space Fractional Partial Differential Equations\thanks{The work of the first author is supported
by National Natural Science Foundation of China (No.
10901027). The work of second author is supported by the National Basic Research Program under the Grant
2011CB30971 and National Natural Science Foundation of China (No. 11171335, 11225107).}}

\author{Yingjun Jiang\thanks{Department of Mathematics and Scientific
Computing, Changsha University of Science and Technology, Changsha,
410076, China (jiangyingjun@csust.edu.cn).} and Xuejun
Xu\thanks{Institute of Computational Mathematics and
Scientific/Engineering Computing, Academy of Mathematics and Systems
Science, Chinese Academy of Sciences, P.O. Box 2719, Beijing,
100190, P.R. China ({\tt xxj@lsec.cc.ac.cn}).}}

\date{}

\begin{document}
\maketitle

\begin{abstract}
In this paper, a two-level additive Schwarz preconditioner is proposed for solving the
algebraic systems resulting from the finite element approximations of space
fractional partial differential equations (SFPDEs).
It is shown that the condition number of the preconditioned system is bounded by $C(1+H/\delta)$, where $H$ is the maximum diameter of subdomains and  $\delta$ is the overlap size among the subdomains.
Numerical results are given to support our theoretical findings.

\end{abstract}


\section{Introduction}\label{sec1}
Space fractional partial differential equations have been wildly used to describe the supper-diffusion processes in the natural world (see \cite{Metzler1}). Let $\Omega$ denote a polyhedral
domain in $\mathbb{R}^d$ and $\tilde{M}(z)$ denotes a probability density function on
$S^{d-1}$, where $S^{d-1}=\{z\in
\mathbb{R}^d; ||z||_2=1\}$ and $||\cdot||_2$ denotes the standard  Euclidean norm.
  In this paper, we consider the following multi-dimensional SFPDE(\cite{Meerschaert0})
\begin{eqnarray}\label{eqn source n}
-\int_{S^{d-1}}D^{2\alpha}_{z}u(x)\, \tilde{M}(z)dz+cu(x)=f(x), \quad x\in
\Omega,
\end{eqnarray}
where $1/2< \alpha < 1$, $c\geq 0$ and $D^{2\alpha}_{z}$, which will be given later, denotes
the directional derivative of order $2\alpha$ in the direction $z$. We assume  $\tilde{M}$ is symmetric about origin,
i.e., $\tilde{M}(z)=\tilde{M}(z')$ if $z,z'\in S^{d-1}$ satisfy $z+z'=0$, which
means that the above SFPDE is symmetric.

Actually, the equation (\ref{eqn source n}) is
 an appropriate extension from one dimensional problem
\begin{equation}\label{eqn surce one d}
   -(p\ _{-\infty}D_x^{2\alpha}+q\ _xD_\infty^{2\alpha})u+cu=f,
\end{equation}
  and its corresponding developing equation can be used to describe a general super-diffusion process (see
\cite{Meerschaert0} for details), where $_{-\infty}D_x^{2\alpha},
{_xD_\infty^{2\alpha}}$ denote Riemann-Liouville fractional
derivatives. One special case of (\ref{eqn source n}) is
\begin{equation}\label{eqn dis}
  -\sum\limits_{i=1}^d(p_i\ _{-\infty}D_{x_i}^{2\alpha}+q_i\ _{x_i}D_\infty^{2\alpha})u+cu=f
\end{equation}
with $p_i,q_i\in \mathbb{R}$ satisfying $p_i=q_i$ and $p_i+q_i=1$.
We may find that equation (\ref{eqn dis}) can be obtained from (\ref{eqn source n}) by taking  \begin{equation}\label{eqn add}
\tilde{M}=\sum_{i=1}^d
p_i\delta(z-e_i)+q_i\delta(z+e_i),
\end{equation}
 where $e_i$ is the $i$th
column of identity matrix in $\mathbb{R}^{d\times d}$ and $\delta$
denotes the Dirac function on $S^{d-1}$.

 Extensive numerical methods have already been developed for SFPDEs,
like finite difference methods (see e.g.,
\cite{Beumer1,Meerschaert1,Meerschaert2,Sousa1,Tadjeran1}),
finite element methods (see e.g., \cite{Deng1,Ervin0,Liu1,Ervin2,Roop1}), and
spectral methods \cite{Li1}.  Due to the nonlocal properties of
fractional differential operators, the most important issue for numerical computation of SFPDEs is how to reduce the
computation costs.
Some methods for the reduction have already been designed, like alternating-direction implicit methods
(ADI)\cite{Meerschaert2,Wang2,Wang3}, special iterative methods \cite{Lin1,Pang1,Wang3,Wang4,Wang5,Wang6} and multigrid methods \cite{Zhou1,Pang1}.

The discrete systems $Ax=b$ of SFPDEs usually have the following characteristics: 1). the condition number
of $A$ increases fast, as the mesh becomes fine; 2). the
coefficient matrix $A$ is dense.
As for iterative methods, two issues need to be concerned for efficiency: one is how to construct  good preconditioners for the discrete system $Ax=b$,  which may help us to save the iterative steps; the other is how to reduce the computation cost of each iterative step, for which some papers employ the multiplication of Toeplitz matrices and vectors with $n \log (n)$ computation complexity (see e.g., \cite{Lin1,Pang1,Wang1,Wang3,Wang4,Wang5,Wang6}).
If the iterative step can be carried out in parallel, the efficiency of solving SFPDEs may be significantly improved. The parallelizable algorithms have been wildly used in numerical solutions for PDEs (see e.g. \cite{Toselli1}). Whereas, to the best of our knowledge, no parallelizable algorithms have been designed for SFPDEs.

In this paper, we shall construct a two-level additive preconditioner for the discrete system $Ax=b$ resulting from the finite element approximation of (\ref{eqn source n}), and then use the preconditioned conjugate method (PCG) to solve it.
The preconditioner we construct is almost optimal, i.e., the condition number of the preconditioned system is bounded by $C(1+H/\delta)$, where $H$ is the maximum diameter of subdomains and $\delta$ is the overlap size among the subdomains. Moreover, the preconditioner may be employed in parallel or each step of the PCG can be carried out in parallel.  As a result, the whole numerical solution processes, including the generation of $A,b$ and the multiplication of matrices and vectors, may be conducted in parallel.

 Without loss of generality, we focus on
 the case $d=2$, namely,
 we consider the problem (\ref{eqn source n}) in $\mathbb{R}^2$. For $\Lambda\subset \mathbb{R}^2$,
 denote $L^2(\Lambda)$ the space of all measurable function $v$ on $\Lambda$ satisfying
$\int_\Lambda (v(x))^2dx<\infty,$ and $C^\infty_0(\Lambda)$  the
space of infinitely differentiable functions with compact support in
$\Lambda$. Set
$$(v,w)_\Lambda=\iint_\Lambda vw dxdy,\quad ||v||_\Lambda=(v,v)^{1/2}_\Lambda,$$
and they are abbreviated as $(v,w)$ and $||v||$ respectively if
$\Lambda=\mathbb{R}^2$.

To simplify our statement, we make a convention here: function $v$
defined on a domain $\Lambda\subset\mathbb{R}^2$ also denotes its
extension on $\mathbb{R}^2$ which extends $v$ by zero outside
$\Lambda$. The constant $C$  with or without subscript shall denote
a generic positive constant which may take on different values in
different places. These constants shall always be independent of mesh sizes and  numbers of subdomains.
Following \cite{Xu1},  we also use symbols $\lesssim,\gtrsim $ and
$\approx $ in this paper. That $a_1\lesssim b_1$, $a_2\gtrsim b_2$
and $a_3\approx b_3$ mean that $a_1\leq C_1 b_1$, $a_2\geq C_2 b_2$
and $C_3b_3\leq a_3\leq C'_3 b_3$ for some positive constants $C_1,C_2,C_3$
and $C'_3$.

The rest of the paper is organized as follows: in section 2, the variational problem of (\ref{eqn source n}) and
 its finite element discretization are described; in section \ref{sec3}, the two-level additive preconditioner for the SFPDEs is presented;
in section \ref{sec4}, it is proved that the  preconditioner is almost optimal; finally, in section
5, the numerical results shall be given to support our theoretical
findings.

\section{The model problem and its discretization}
In this section, we shall describe the SFPDEs in details, and
then introduce its variational formulation and finite element
discretization.

\subsection{The model problem} \label{sec}
\begin{definition}\label{defi 1}\cite{Ervin2}
Let $\mu>0$, $\theta\in \mathbb{R}$. The $\mu$th order fractional
integral in the direction $z=(\cos\theta,\sin\theta)$ is defined by
$$D^{-\mu}_z v(x,y):=D^{-\mu}_\theta v(x,y)=\int_0^\infty \frac{\tau^{\mu-1}}{\Gamma(\mu)}v(x-\tau \cos\theta,y-\tau\sin\theta)d\tau,$$
where $\Gamma$ is the Gamma function.
\end{definition}

\begin{definition}\label{defi 2}\cite{Ervin2}
Let $n$ be a positive integer, and $\theta\in \mathbb{R}$. The $n$th
order derivative in the direction of $z=(\cos\theta,\sin\theta)$ is
given by
$$D^n_\theta v(x,y):=\left(\cos\theta\frac{\partial}{\partial x}+\sin\theta\frac{\partial}{\partial y}\right)^nv(x,y).$$
\end{definition}
\begin{definition}\label{defi 3}\cite{Ervin2}
Let $\mu>0$, $\theta\in \mathbb{R}$. Let $n$ be the integer such
that $n-1\leq\mu<n$, and define $\sigma=n-\mu$. Then the $\mu$th
order directional derivative in the direction of
$z=(\cos\theta,\sin\theta)$ is defined by
$$D^\mu_z v(x,y):=D^\mu_\theta v(x,y)=D^n_\theta D_\theta^{-\sigma}v(x,y).$$
\end{definition}
If  $v$ is viewed as a function in $x$,
$D^{\mu}_0$, $D^{\mu}_\pi$ are
just the left and the right Riemamm-Liouville derivatives (see Appendix).
\begin{definition}\label{defi 4}\cite{Ervin2}
Assume that $v:\mathbb{R}^2\rightarrow \mathbb{R}$, $\mu>0$. The $\mu$th
order fractional derivative with respect to the measure $\tilde{M}$ is
defined as
$$D^\mu_{\tilde{M}}v(x,y):=\int_{S^1}D^\mu_\theta v(x,y)\tilde{M}(\theta)d\theta,$$
where $S^1=[0+\nu, 2\pi+\nu)$ with a suitable scalar $\nu$, and
$\tilde{M}(\theta)$, which satisfies
$\int_{\nu}^{2\pi+\nu}\tilde{M}(\theta)d\theta=1$, is a periodic function
with period $2\pi$. Usually we take $\nu=0$, if it causes no
unreasonable expression (see (\ref{eqn add})).
\end{definition}

For $u:\mathbb{R}^2\rightarrow \mathbb{R}$, define differential
operator $L_\alpha$ in $\mathbb{R}^2$ as
$$L_\alpha u=-D^{2\alpha}_{\tilde{M}} u+cu.$$
Denote $\Omega$ a polygonal domain in $\mathbb{R}^2$, set
$1/2<\alpha< 1$. The model problem of this paper
is to find $u:\bar\Omega\rightarrow \mathbb{R}$ such that
\begin{equation}\label{eqn source new}
\left\{\begin{array}{cc}L_\alpha u=f, & \hbox{in}\ \Omega,\\
 u=0, &  \hbox{on}\
\partial\Omega,\end{array}\right.
\end{equation}
where $f$ is a source term and we assume that $\tilde{M}(\theta)$ satisfies
$\tilde{M}(\theta)=\tilde{M}(\theta+\pi)$ for $\theta\in \mathbb{R}$, i.e.,
(\ref{eqn source new}) is a symmetric problem. Here, we recall the
convection made in Section \ref{sec1}:   $u$ also denotes its
extension by zero outside $\Omega$.

\subsection{The variational formulation and finite element discretization}
\begin{definition}\label{defi 8}\cite{Tartar1} Let $\mu\geq 0$, $\mathcal{F}{v}(\xi_1,\xi_2)$ be the Fourier transform of $v(x,y)$, $|\xi|=\sqrt{\xi_1^2+\xi_2^2}$.
Define norm
$$||v||_{H^\mu(\mathbb{R}^2)}:=\left\|(1+|\xi|^2)^{\mu/2}|\mathcal{F}{v}|\right\|.$$
Let $H^\mu(\mathbb{R}^2):=\{v\in L^2(\mathbb{R}^2);
||v||_{H^\mu(\mathbb{R}^2)}< \infty\}$.
\end{definition}
For $v\in H_0^\mu(\Omega)$,  we also denote  $||v||_{H^\mu(\mathbb{R}^2)}$ by $||v||_{H^\mu(\Omega)}$. It is known that $H^\mu(\mathbb{R}^2)$ is a Hilbert space equipped
with the inner product
$(v,w)_{H^\mu(\mathbb{R}^2)}=((1+|\xi|^2)^{\mu}\mathcal{F}{v},\overline{\mathcal{F}{w}})$
and $C_0^\infty(\mathbb{R}^2)$ is dense in $H^\mu(\mathbb{R}^2)$
(see \cite{Tartar1}).
Since we employ the finite element discretization, the weak fractional directional
derivative need to be introduced. Let
$L_{loc}^1(\mathbb{R}^2)$ denote the set of locally integrable
functions on $\mathbb{R}^2$.

\begin{definition}\cite{Jiang1}\label{defi 0} Given $\mu>0$,
$\theta\in \mathbb{R}$, let $v\in L^2(\mathbb{R}^2)$. If there is a
function $v_\mu\in L^1_{loc}(\mathbb{R}^2)$ such that
$$(v,{D^\mu_{\theta+\pi}}w)=({v_\mu},w),\quad \forall w\in C_0^\infty(\mathbb{R}^2),$$
then $v_\mu$ is called the weak $\mu$th order derivative in the
direction of $\theta$ for $v$, denoted by ${D^\mu_{\theta}}v$, i.e.,
$v_\mu={D^\mu_{\theta}}v$.
\end{definition}

The weak derivative ${D^\mu_{\theta}}v$
is unique if it exists and the weak derivative coincides with
the correspondent derivative defined in Definition \ref{defi 3} if $v\in
C^\infty_0(\mathbb{R}^2)$.  In the following, we use
${D^\mu_{\theta}}v$ to denote the weak derivative.

\begin{lemma}\cite{Ervin2,Jiang1}\label{lemma a9}
Let $\mu>0$. For any $v\in H^\mu(\mathbb{R}^2)$, $0<s\leq\mu$ and $\theta\in \mathbb{R}$, the weak derivative
${D^s_{\theta}}v$ exists and satisfies
\begin{equation}\label{r1}\mathcal{F}D^s_\theta v(\xi_1,\xi_2)=(2\pi i\xi_1\cos \theta+2\pi i\xi_2\sin\theta)^s \mathcal{F}v(\xi_1,\xi_2),\end{equation}
\begin{equation}\label{r2}||D^s_\theta v||\leq C||v||_{H^\mu(\mathbb{R}^2)}.\end{equation}
\end{lemma}

Define the bilinear form $\tilde{B}: H^\alpha_0(\Omega)\times
H^\alpha_0(\Omega)\rightarrow \mathbb{R}$ as
$$
\tilde{B}(u,v):=-\int_0^{2\pi}(D^\alpha_\theta u,D^\alpha_{\theta+\pi}
v)\tilde{M}(\theta)d\theta+c(u,v).
$$
Because $\tilde{M}(\theta)=\tilde{M}(\theta+\pi)$ for $\theta\in \mathbb{R}$, it is easy
to check that $\tilde{B}(v,w)$ is a symmetric bilinear  form, i.e.,
$\tilde{B}(v,w)=\tilde{B}(w,v)$ for $v,w\in H^{\alpha}_0(\Omega)$.
The variational formulation of (\ref{eqn source new}) (see
\cite{Ervin2,Jiang1}) is to find $u\in H^{\alpha}_0(\Omega)$ such that
\begin{equation}\label{eqn variational form new}
\tilde{B}(u,v)=(f,v),\quad \forall v\in H^{\alpha}_0(\Omega).
\end{equation}
Here $\tilde{M}$ is taken such that (\ref{eqn variational form new}) admits a unique solution in $H^{\alpha}_0(\Omega)$ (for the details, we refer to \cite{Ervin2,Jiang1}).

We construct a quasi-uniform triangulation $\Gamma_H=\{\Omega_i\}_{i=1}^J$ of $\Omega$ with
$$\hbox{the diameter of } \Omega_i \approx O(H).$$
Divide each $\Omega_i$ into smaller simplices $\tau_j$ of diameter $O(h)$, such that $\Gamma_h=\{\tau_j\}$ form a finer triangulation of $\Omega$.
Denote $V_H$ and $V_h$ piecewise linear finite element function spaces defined on the triangulations $\Gamma_H$ and $\Gamma_h$ respectively.
 It is known that
$V_h\subset H^1_0(\Omega)\subset
H^\alpha_0(\Omega)$.

The finite element approximation for (\ref{eqn
variational form new}) (the details please see  \cite{Jiang1}) is to find $u_h\in V_h$ such that
 \begin{equation}\label{eqn finite element new}
 B(u_h,v)=(f,v),\quad \forall v\in V_h,
 \end{equation}
where $B(v,w)=-\int_0^{2\pi}(D^\alpha_\theta v,D^\alpha_{\theta+\pi}
w)M(\theta)d\theta+c(v,w)$, $M(\theta)$ is equal to a discrete form $\sum_{k=1}^Lp_k\delta(\theta-\theta_k)$ such that $B(\cdot,\cdot)$ is a symmetric
bilinear form, and
\begin{equation}\label{r5r6}
B(v,v)\gtrsim ||v||^2_{H^{\alpha}(\Omega)}, \quad B(v,w)\lesssim ||v||_{H^{\alpha}(\Omega)}||w||_{H^{\alpha}(\Omega)},\quad v,w\in H^{\alpha}_0(\Omega).
\end{equation}
Meanwhile
 $$\int_0^{2\pi}M(\theta)d\theta\lesssim1.$$

\begin{remark}\label{remark 1}
The direct finite element discretization of (\ref{eqn
variational form new}) is
 \begin{equation}\label{finite}
\tilde{B}(u_h,v)=(f,v),\quad v\in V_h.
\end{equation}
But the discretization is hardly computed. So we use (\ref{eqn finite element new}) instead of (\ref{finite}), where $B(\cdot,\cdot)$ is understood as the approximation to $\tilde{B}(\cdot,\cdot)$. For the details, please refer to \cite{Jiang1}.
\end{remark}
\section{A two-level additive Schwarz preconditioner}\label{sec3}
Take $f_h\in V_h$ such that $(f_h,v)=(f,v)$, $\forall v\in V_h$ and
define a linear operator $A:V_h\rightarrow V_h$ satisfying
\begin{equation}\label{eqn variational form 2}
(Av,w)=B(v,w), \quad \forall v,w\in V_h.
\end{equation}
Since $B(v,w)$ is a symmetric
bilinear form, by (\ref{r5r6}), we know that
$A:V_h\rightarrow V_h$ is symmetric positive definite with respect to
$(\cdot,\cdot)$, i.e.,
$$(Av,w)=(v,Aw),\quad v,w\in V_h; \quad (Av,v)>0,\quad 0\neq v\in V_h.$$
Then bilinear form
$$(v,w)_A:=(Av,w)$$
also induces an inner product on $V_h$. Set norm
$$||v||_A=(Av,v)^{1/2},\quad v\in V_h.$$
By (\ref{r5r6}), we have
\begin{equation}\label{eqn aaaa5}
||v||_A\approx ||v||_{H^\alpha(\Omega)}, \quad \forall v\in V_h.
\end{equation}
The problem (\ref{eqn finite element new}) can be
restated as to find $u_h\in V_h$ such that
\begin{equation}\label{M eqn problem111}
Au_h=f_h.
\end{equation}
For the above equation, we shall construct our two-level Schwarz preconditioner and then use PCG method to solve it.

Our preconditioner is designed by making use of the following overlapping domain decomposition $\Omega=\cup_{i=1}^J\{\Omega'_i\}$, where the subdomain $\Omega'_i$ contains coarse subdomain $\Omega_i$, and satisfies
 $$\hbox{the diameter of } \Omega'_i \approx O(H).$$
Meanwhile the boundary of $\Omega'_i$ align with the mesh of triangulation $\Gamma_h$, and
 the distance from $\partial \Omega'_i\cap \Omega$ to  $\Omega_i$ is greater than $\delta$ , which is a positive constant measuring the  overlapping size among the subdomians.
 Define subspaces $V_i$ of $V_h$ as
 $$V_i=\{v\in V_h| v(x)=0,x\in \Omega\setminus\Omega'_i\},\quad i=1,2,\ldots,J,$$
 and let $V_0=V_H$.

For our analysis, we regroup the subregions in terms of the coloring strategy (see e.g. \cite{Toselli1}). By a minimal or good coloring, we group  the subregions $\{\Omega'_i\}$ into $J^C$ classes, each of which has some disjoint subregions and can be regarded as one subregion. Exactly, decompose the index set $\{1,2,\ldots,J\}=\cup_{i=1}^{J^C}I_i$ with $I_i$ satisfying that $\Omega'_l\cap \Omega'_k=\emptyset$ for any $l, k\in I_i (l\neq k)$; for $i=1,\ldots,J^C$, define new subregions $\tilde{\Omega}_i=\cup_{j\in I_i} \Omega'_j$ and new subspaces $\tilde{V}_i=\bigoplus_{j\in I_i} V_j$.

For each $k\in \{0,1,2,\ldots,J\}$, we define some projectors
$Q_k,P_k:V_h\rightarrow V_k$  by
$$(Q_kv,w)=(v,w),\quad (P_kv,w)_A=(v,w)_A,\quad \forall w\in V_k,v\in V_h,$$
 and define the
linear operator $A_k:V_k\rightarrow V_k$ by
$$(A_kv,w)=(Av,w),\quad v,w\in V_k.$$
It is not hard to verify that
\begin{equation}\label{M eqn important}
A_kP_k=Q_kA,\quad k=0,1,\ldots,J.
\end{equation}
To help our analysis, we define some projectors $\tilde{P}_i:V_h\rightarrow \tilde{V}_i$, $i=1,\ldots,J^C$, by
$$(\tilde{P}_iv,w)_A=(v,w)_A,\quad \forall w\in \tilde{V}_i, v\in V_h.$$
\begin{remark}\label{remark 11}
Different from the integer order PDEs, it is interesting to see that $\tilde{P}_i\neq\sum_{j\in I_i}P_j.$
\end{remark}

Now, we are ready to present our two-level additive Schwarz pre-conditioner, i.e.,
$$B_h=\sum_{i=0}^J A_i^{-1}Q_i.$$
 By (\ref{M eqn important}), we have
$$B_hA=\sum_{i=0}^J A_i^{-1}Q_iA=\sum_{i=0}^J A_i^{-1}A_iP_i=\sum_{i=0}^J P_i.$$
Define $P_h: =\sum_{i=0}^J P_i$, and then the  preconditioned system is
\begin{equation}\label{preconditioner system}
P_hu_h=B_hf_h.
\end{equation}

In the next section, we shall prove the condition number of $P_h$ is bounded by $C(1+H/\delta)$,
where the constant $C$ is independent of mesh size and the numbers of subdomains, but dependent of $J^C$.

\section{ Condition number estimate}\label{sec4}

We first introduce two interpolation norms and relevant Sobolev spaces (see e.g., \cite{Tartar1}).
Let $\Lambda$ be a
domain in $\mathbb{R}^2$. For integer $m$, denote by $||\cdot||_{\tilde{H}^{m}(\Lambda)}$ the Sobolev norm of integer order $m$, i.e.,
$$||v||_{\tilde{H}^m(\Lambda)}:=\left(\sum\limits_{|l|\leq m}||D^lv||_{L^2(\Lambda)}^2\right)^{1/2},$$
with $l=(l_1,l_2)$, $|l|=l_1+l_2$ and
$D^l=(\frac{d}{dx})^{l_1}(\frac{d}{dy})^{l_2}$. Let $\mu>0$ be a non-integer and $0<s<1$,  $n$ is a non-negative integer such that $n<\mu<n+1$. We introduce the interpolation norms
\begin{equation}\label{norm 1}||v||_{\tilde{H}^{\mu}(\Lambda)}:=\left(\int_0^{\infty}\tilde{K}(v,t) t^{-2\mu-1}dt\right)^{1/2},||v||_{\hat{H}^{s}(\Lambda)}:=\left(\int_0^{\infty}\hat{K}(v,t) t^{-2s-1}dt\right)^{1/2},
\end{equation}
where
$$\tilde{K}(v,t):=\inf_{w\in \tilde{H}^{n+1}(\Lambda)}\left(||v-w||^2_{\tilde{H}^n(\Lambda)}+t^2||w||^2_{\tilde{H}^{n+1}(\Lambda)}\right),$$
$$\hat{K}(v,t):=\inf_{w\in \tilde{H}_0^{1}(\Lambda)}\left(||v-w||^2_{L^2(\Lambda)}+t^2||w||^2_{\tilde{H}^{1}(\Lambda)}\right).$$
 Relevant Sobolev spaces are
\begin{equation}\label{eqn aaaaaa3}
{\tilde{H}}^{\mu}(\Lambda):=\{v\in
L^2(\Lambda);||v||_{{\tilde{H}}^{\mu}(\Lambda)} <\infty\},\quad{\hat H}^{s}(\Lambda):=\{v\in
L^2(\Lambda);||v||_{{\hat{H}}^{s}(\Lambda)} <\infty\}.
\end{equation}

For $\mu>0$, it is known that $\tilde{H}^\mu(\mathbb{R}^2)$ coincides with
${H}^\mu(\mathbb{R}^2)$.  The following norms relation is
useful in our analysis.

\begin{lemma}\label{lemma aaaab} \cite{Lions1,Tartar1} Let $0<\mu<1$ with $\mu\neq 1/2$ and $\Lambda$ be a domain in $\mathbb{R}^2$ with Lipschitz  boundary. Then $\hat{{H}}^\mu(\Lambda)$ coincides with $\tilde{H}^\mu_0(\Lambda)$ with equivalent norms.
\end{lemma}

In the following, we shall give some useful results.

\begin{lemma}\label{lemma d1}
If  $\hbox{dist}(\Omega'_i,\Omega'_j)\geq lH$ for integer $l\geq1$,  we have
\begin{equation}\label{eqn d1}(D^{\alpha}_\theta v,D^{\alpha}_{\theta+\pi}
w)\lesssim \frac{1}{l^{0.5+\alpha}}||v||_A ||w||_A, \quad  v\in V_i, w\in V_j.
\end{equation}
\end{lemma}
\begin{figure}[h]
\begin{center}
  \includegraphics[width=10cm]{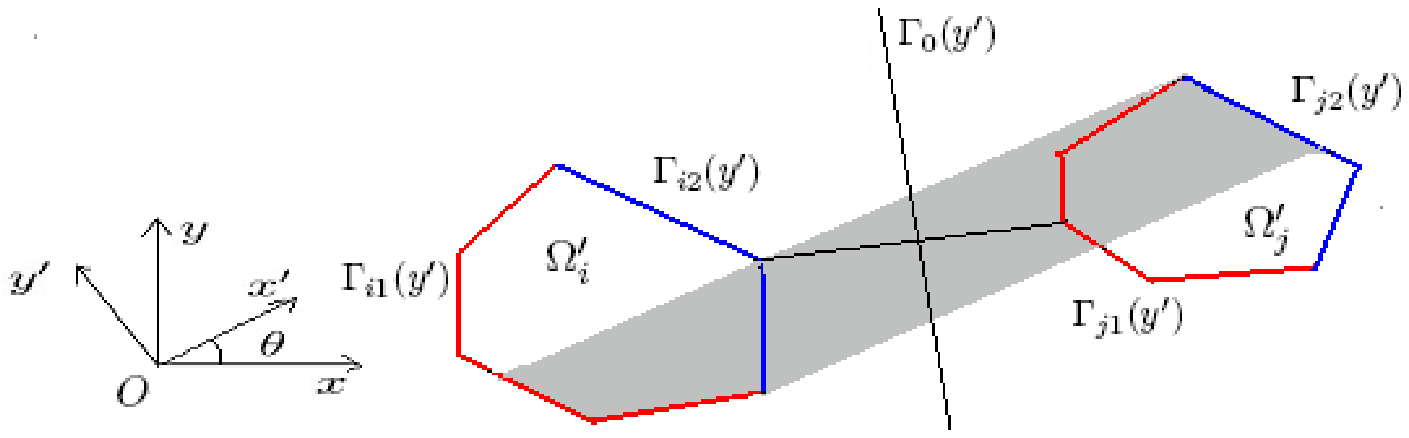}\\
  \caption{}\label{fig 1}
  \end{center}
\end{figure}

\proof Without loss of generality, we prove that (\ref{eqn d1}) holds under the situation as the figure \ref{fig 1} shows. Cartesian coordinate $x'Oy'$ is obtained by rotating $xOy$ $\theta$ angle counterclockwise. We set that $p_i\in \bar\Omega_i$ and $p_j\in \bar\Omega_j$ such that the length of line segment
 $p_ip_j$ is equal to $\hbox{dist}(\Omega'_i,\Omega'_j)$, the graphs of functions $x'=\Gamma_{i1}(y')$ and $x'=\Gamma_{i2}(y')$ are parts of boundary $\partial \Omega'_i$ such that $\Omega'_i=\{(x',y'):\Gamma_{i1}(y')<x'<\Gamma_{i2}(y')\}$,  the graphs of functions $\Gamma_{j1}(y')$ and $\Gamma_{j2}(y')$ are parts of boundary $\partial \Omega'_j$ such that $\Omega'_j=\{(x',y'):\Gamma_{j1}(y')<x'<\Gamma_{j2}(y')\}$,  and the graph of $x'=\Gamma_0(y')$ is the perpendicular bisector of $p_ip_j$.  Then we have $D^{\alpha}_\theta={_{-\infty}}D^{\alpha}_{x'}$, $D^{\alpha}_{\theta+\pi}={_{x'}D^{\alpha}_{\infty}}$, where ${_{-\infty}}D^{\alpha}_{x'}$ and ${_{x'}D^{\alpha}_{\infty}}$ are the left and right Riemann-Liouville fractional derivative operator (see Appendix).
 We in the following only prove (\ref{eqn d1}) for the case that there exists a ray in $(\cos\theta,\sin\theta)$ direction going through both $\Omega'_i$ and $\Omega'_j$ from $\Omega'_i$ to $\Omega'_j$. Indeed, when the other case is true, $(D^{\alpha}_\theta v,D^{\alpha}_{\theta+\pi}w)=0$ and (\ref{eqn d1}) naturally holds.

Denote
$$M_{y'}=\max\limits_{\tiny{\begin{array}{c}(x',y')\in \Omega'_i\\
(x'',y')\in \Omega'_j\end{array}}}y',\quad m_{y'}=\min\limits_{\tiny{\begin{array}{c}(x',y')\in \Omega'_i\\
(x'',y')\in \Omega'_j\end{array}}}y',$$
$$M_{x'}=\max\limits_{\tiny{\begin{array}{c}
(x',y')\in \Omega'_j\end{array}}}x',\quad m_{x'}=\min\limits_{\tiny{\begin{array}{c}(x',y')\in \Omega'_i\end{array}}}x'.$$
The shaded area is
$$\Lambda=\{(x',y')| \Gamma_{i,1}(y')\leq x'\leq \Gamma_{j,2}(y'),m_{y'}\leq y'\leq M_{y'}\}.$$
 Let $\Lambda_i=\{(x',y')|x'\leq \Gamma_0(y')\}\cap \Lambda$  and $\Lambda_j=\{(x',y')|x'\geq \Gamma_0(y')\}\cap \Lambda$. Then we have
\begin{eqnarray}\label{eqn use}
(D^{\alpha}_\theta v,D^{\alpha}_{\theta+\pi}w)&=&(D^{\alpha}_\theta v,D^{\alpha}_{\theta+\pi}w)_\Lambda=(D^{\alpha}_\theta v,D^{\alpha}_{\theta+\pi}w)_{\Lambda_i}+(D^{\alpha}_\theta v,D^{\alpha}_{\theta+\pi}w)_{\Lambda_j}\nonumber\\
&\leq& ||D^{\alpha}_{\theta}v||_{\Lambda_i} ||{D^{\alpha}_{\theta+\pi}}w||_{\Lambda_i}+||D^{\alpha}_{\theta}v||_{\Lambda_j} ||{D^{\alpha}_{\theta+\pi}}w||_{\Lambda_j}.
\end{eqnarray}
For $(x,y)\in \Lambda_j$ whose coordinate under $x'Oy'$ is $(x',y')$, noting that $\hbox{supp}(v)\subset \Omega'_i$, we have
\begin{eqnarray}\label{eqn d2}
D^{\alpha}_{\theta}v(x,y)&=&{_{-\infty}}D^{\alpha}_{x'}v(x',y')
=\frac{d}{dx'}{_{-\infty}}D^{-(1-\alpha)}_{x'}v(x',y')\nonumber\\
&=&\frac1{\Gamma(1-\alpha)}\frac{d}{dx'}\int_{-\infty}^{x'}(x'-s)^{-\alpha}v(s,y')ds\nonumber\\
&=&\frac1{\Gamma(1-\alpha)}\frac{d}{dx'}\int_{\Gamma_{i1}(y')}^{\Gamma_{i2}(y')}(x'-s)^{-\alpha}v(s,y')ds\nonumber\\
&=&\frac1{\Gamma(-\alpha)}\int_{\Gamma_{i1}(y')}^{\Gamma_{i2}(y')}(x'-s)^{-(1+\alpha)}v(s,y')ds\nonumber\\
&\lesssim&(lH)^{-(1+\alpha)}\int_{\Gamma_{i1}(y')}^{\Gamma_{i2}(y')}|v(s,y')|ds\nonumber\\
&\lesssim&(lH)^{-(1+\alpha)}H^{1/2}\left(\int_{\Gamma_{i1}(y')}^{\Gamma_{i2}(y')}v^2(s,y')ds\right)^{1/2},
\end{eqnarray}
where in the fifth equality we have used the relation $\Gamma(1-\alpha)=-\alpha\Gamma(-\alpha)$, the first inequality is by $(x'-s)\geq lH/2$ when $s\leq \Gamma_{i2}(y')$ and the last by the Cauchy-Schwarz inequality.
Then
\begin{eqnarray}\label{eqn d3}
||D^{\alpha}_{\theta}v||^2_{\Lambda_j}&\lesssim& (lH)^{-2(1+\alpha)}H\int_{m_{x'}}^{M_{x'}}dx'\int_{m_{y'}}^{M_{y'}}dy'\int_{\Gamma_{i1}(y')}^{\Gamma_{i2}(y')}v^2(s,y')ds\nonumber\\
&\lesssim& (lH)^{-(1+2\alpha)}H\int_{m_{y'}}^{M_{y'}}dy'\int_{\Gamma_{i1}(y')}^{\Gamma_{i2}(y')}v^2(s,y')ds\nonumber\\
&\leq&(lH)^{-(1+2\alpha)}H\iint_{\Omega'_i}v^2dx'dy'=l^{-(1+2\alpha)}H^{-2\alpha}||v||^2_{\Omega'_i},
\end{eqnarray}
where the second inequality is by $\int_{m_{x'}}^{M_{x'}}dx'\lesssim lH$.
Denote $d_i$ as the diameter of $\Omega'_i$, and define a function in $x'$ as
$$H_i(x')=\left\{\begin{array}{ll} 0,& \hbox{if } x'> d_i;\\
\frac{x'^{(\alpha-1)}}{\Gamma(\alpha)} &\hbox{if }0\leq x'\leq d_i.
\end{array}\right.$$
For $(x,y)\in \Omega'_i$ (whose coordinate is $(x',y')$ under $x'Oy'$),  by (\ref{property 1}), we have
\begin{eqnarray}\label{eqn d5}
v(x,y)=D^{-\alpha}_{\theta}D^{\alpha}_{\theta}v(x,y)={_{-\infty}}D_{x'}^{-\alpha}\ {_{-\infty}}D_{x'}^{\alpha}v(x',y')=H_i*{_{-\infty}}D_{x'}^{\alpha}v(\cdot,y'),
\end{eqnarray}
where $v*w$ denote the convolution product (see e.g., \cite{Adams1}). Then by the Young Theorem (Theorem 4.30 in \cite{Adams1}),
\begin{equation}\label{eqn d6}
||v(\cdot,y')||_{[\Gamma_{i1},\Gamma_{i2}]}\leq ||v(\cdot,y')||_{\mathbb{R}}\leq ||H_i||_{L^1(\mathbb{R})}||{_{-\infty}}D_{x'}^{\alpha}v(\cdot,y')||_{\mathbb{R}}\lesssim H^\alpha||{_{-\infty}}D_{x'}^{\alpha}v(\cdot,y')||_{\mathbb{R}}.
\end{equation}
Furthermore we have
\begin{eqnarray}\label{eqn d4}
||v||^2_{\Omega'_i}&=&\iint_{\Omega'_i}v^2dx'dy'=\int^{M_{y'}}_{m_{y'}}||v(\cdot,y')||^2_{[\Gamma_{i1}(y'),\Gamma_{i2}(y')]}dy'\nonumber\\
&\lesssim& H^{2\alpha} ||{_{-\infty}}D_{x'}^{\alpha}v||^2_{\mathbb{R}^2}=H^{2\alpha}||D_{\theta}^{\alpha}v||^2_{\mathbb{R}^2}\nonumber\\
&\lesssim&H^{2\alpha}||v||^2_{H^\alpha(\Omega)}\approx H^{2\alpha}||v||^2_A,
\end{eqnarray}
where the second inequality is by Lemma \ref{lemma a9}, the last equality is by (\ref{eqn aaaa5}).
Combining (\ref{eqn d4}) with (\ref{eqn d3}), we obtain
\begin{equation}\label{eqn d7}
||D^{\alpha}_{\theta}v||_{\Lambda_j}\lesssim \frac1{l^{1/2+\alpha}}||v||_A.
\end{equation}
Similarly we may also obtain
\begin{equation}\label{eqn d8}
||D^{\alpha}_{\theta+\pi}w||_{\Lambda_i}\lesssim \frac1{l^{1/2+\alpha}}||w||_A.
\end{equation}
By Lemma \ref{lemma a9},
\begin{equation}\label{eqn d9}
||D^{\alpha}_{\theta+\pi}w||_{\Lambda_j}\lesssim ||w||_{H^\alpha(\Omega)}\approx ||w||_A,\quad ||D^{\alpha}_{\theta}v||_{\Lambda_i}\lesssim ||v||_{H^\alpha(\Omega)}\approx ||v||_A.
\end{equation}
Combining with (\ref{eqn use}), (\ref{eqn d7}) and (\ref{eqn d8}), we may obtain (\ref{eqn d1}). \eproof

\begin{lemma}\label{lemma d2} Let $S\subset \{1,2,\ldots,J\}$ denote an index set such that $\Omega'_i\cap \Omega'_j =\emptyset$ for any $i,j\in S$, $i\neq j$.
Then for $v=\sum_{i\in S}v_i$ with $v_i\in V_i$, we have
\begin{equation}\label{eqn dd1}
(v,v)_A\lesssim \sum\limits_{i\in S} (v_i,v_i)_A.
\end{equation}
\end{lemma}
\proof
Note that
$$(v,v)_A=B(v,v)=-\int_0^{2\pi}(D^\alpha_\theta v,D^\alpha_{\theta+\pi}
v)M(\theta)d\theta+c(u,v),$$  it is easy to see  that the lemma follows after we prove that
\begin{equation}\label{eqn dd2}
|(D^\alpha_\theta v,D^\alpha_{\theta+\pi}v)|\leq C\sum\limits_{i\in S} (v_i,v_i)_A
\end{equation}
holds for any $\theta$. So next, we give a proof of (\ref{eqn dd2}).
It is easy to see that
\begin{eqnarray}\label{eqn dd3}
|(D^\alpha_\theta v,D^\alpha_{\theta+\pi}v)|&\leq& \sum\limits_{k,j\in S} |(D^\alpha_\theta v_k,D^\alpha_{\theta+\pi}v_j)|.
\end{eqnarray}
In fact,  we may prove
\begin{equation}\label{eqn dd4}
\sum\limits_{k,j\in S} |(D^\alpha_\theta v_k,D^\alpha_{\theta+\pi}v_j)|\leq C\sum\limits_{i\in S}(v_i,v_i)_A
\end{equation}
 through using the following two kind of inequalities:
\begin{eqnarray}\label{eqn dd9}
 (D^\alpha_\theta v_k,D^\alpha_{\theta+\pi}v_j)&\leq& ||D^\alpha_\theta v_k||\, ||D^\alpha_{\theta+\pi}v_j||
 \nonumber\\
 &\leq& C ||v_k||_A||v_j||_A
\leq  C ||v_k||^2_A+C||v_j||^2_A,\hbox{ if } dist(\Omega'_k,\Omega'_j)< H
\end{eqnarray}
   (by Lemma \ref{lemma a9} and  (\ref{eqn aaaa5}));
\begin{equation}\label{eqn ddd1}
(D^{\alpha}_\theta v_k,D^{\alpha}_{\theta+\pi}
v_j)\leq \frac{C}{l^{0.5+\alpha}}||v_k||_A ||v_j||_A\leq \frac{C}{l^{0.5+\alpha}}(||v_k||^2_A+||v_j||^2_A), \hbox{ if }dist(\Omega'_k,\Omega'_j)\geq lH
 \end{equation} (by Lemma \ref{lemma d1}).

 Fixing $i\in S$, we write
\begin{equation}\label{eqn dd5}
S_i=\{j\in S:\ dist(\Omega'_i,\Omega'_j)< H\},\quad S'_i=\{j\in S:\ dist(\Omega'_i,\Omega'_j)\geq H\}.
\end{equation}
It is easy to see that
\begin{eqnarray}\label{eqn dd6}
&&\hbox{the sum of the terms containing $v_i$ on the left hand of (\ref{eqn dd4})}\nonumber\\
&\leq&\sum\limits_{j\in S} |(D^\alpha_\theta v_i,D^\alpha_{\theta+\pi}v_j)|+\sum\limits_{j\in S} |(D^\alpha_\theta v_j,D^\alpha_{\theta+\pi}v_i)|\nonumber\\
&=&\sum\limits_{j\in S_i} |(D^\alpha_\theta v_i,D^\alpha_{\theta+\pi}v_j)|+\sum\limits_{j\in S_i} |(D^\alpha_\theta v_j,D^\alpha_{\theta+\pi}v_i)|\nonumber\\
&&\quad +\sum\limits_{j\in S'_i} |(D^\alpha_\theta v_i,D^\alpha_{\theta+\pi}v_j)|+\sum\limits_{j\in S'_i} |(D^\alpha_\theta v_j,D^\alpha_{\theta+\pi}v_i)|.
\end{eqnarray}
We know that
\begin{equation}\label{eqn cc1}
card(S_i)\leq C,
\end{equation}
where $card(S)$ denotes the number of elements contained in $S$, $C$ is a positive constant dependent of $J^C$.
By (\ref{eqn cc1}) and (\ref{eqn dd9}), we have
\begin{equation}\label{eqn dd8}
\sum\limits_{j\in S_i} |(D^\alpha_\theta v_i,D^\alpha_{\theta+\pi}v_j)|+\sum\limits_{j\in S_i} |(D^\alpha_\theta v_j,D^\alpha_{\theta+\pi}v_i)|\leq C ||v_i||^2_A+\hbox{terms which do not contain } v_i.
\end{equation}
\begin{figure}[h]
\begin{center}
  \includegraphics[width=9cm]{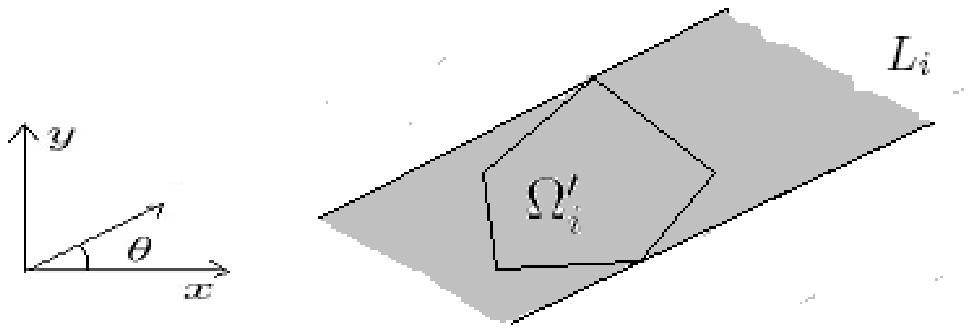}\\
  \caption{}\label{fig 2}
\end{center}
\end{figure}

 As figure \ref{fig 2} shows, $L_i$, containing $\Omega'_i$, denotes the shaded bar area whose edges  are parallel to the vector $(\cos(\theta),\sin(\theta))$ and go through $\bar\Omega'_i$.
It is not hard to see that $(D^\alpha_\theta v_i,D^\alpha_{\theta+\pi}v_j)=0$ if $\Omega'_j\bigcap L_i=\emptyset$, and for positive integer $l$,  the number of the elements of the set $S'_{i,l}=\{\Omega'_k; \Omega'_k\cap L_i\neq\emptyset, lH\leq\hbox{dist}(\Omega'_k, \Omega'_i)< (l+1)H \}$  is bounded by a positive $C$  which is a positive constant dependent of $J^C$. And then using (\ref{eqn ddd1}), we have
\begin{eqnarray}\label{eqn ddd2}
&&\sum\limits_{j\in S'_i} |(D^\alpha_\theta v_i,D^\alpha_{\theta+\pi}v_j)|+\sum\limits_{j\in S'_i} |(D^\alpha_\theta v_j,D^\alpha_{\theta+\pi}v_i)|\nonumber\\
&\leq& \sum\limits_{l\geq 1}\sum\limits_{j\in S'_{i,l}} |(D^\alpha_\theta v_i,D^\alpha_{\theta+\pi}v_j)|+\sum\limits_{l\geq 1}\sum\limits_{j\in S'_{i,l}} |(D^\alpha_\theta v_j,D^\alpha_{\theta+\pi}v_i)|\nonumber\\
&\leq& \sum\limits_{l=1}^{\infty}\frac{C}{l^{0.5+\alpha}}||v_i||^2_A+\hbox{terms which do not contain } v_i\nonumber\\
&\leq& {C}||v_i||^2_A+\hbox{terms which do not contain } v_i.
\end{eqnarray}
From above analysis, using  (\ref{eqn dd9}) and (\ref{eqn ddd1}) to enlarge the correspondent terms of the left hand of (\ref{eqn dd4}), we may show that the inequality (\ref{eqn dd4}) holds with constant $C$ being independent of $\theta$.\eproof

\begin{lemma}\label{lemma d3}
There exist linear operators $\kappa_i$ ($i=0,1,\ldots,J$), $\kappa_0$ maps $\tilde{H}_0^1(\Omega)$ into $V_0$, $\kappa_i(i=1,\ldots,J)$ maps $\tilde{H}_0^1(\Omega)$ and $V_h$ into $\tilde{H}_0^1(\Omega'_i)$ and $V_i$ respectively, such that
   for any $v\in {\tilde{H}_0^1(\Omega)}$,
\begin{equation}\label{eqn dddd3}
v=\sum_{i=0}^J \kappa_i v,
\end{equation}
\begin{equation}\label{eqn ddd10}
\sum\limits_{i=0}^J||\kappa_i v||^2\leq  C||v||^2
\end{equation}
and
\begin{equation}\label{eqn ddd9}
\sum\limits_{i=0}^J||\kappa_i v||^2_{\tilde{H}^1(\Omega)}\leq C\left(1+\frac{H}{\delta}\right)||v||^2_{\tilde{H}^1(\Omega)},
\end{equation}
where the constant $C$ depends on $J^C$.
\end{lemma}
\proof We employ the symbols used in the Lemma 3.12 of \cite{Toselli1}. First let $\tilde I^H:H^1_0(\Omega)\rightarrow V_H$  be a quasi-interpolation operator, which is defined as, for the node $y$ in the coarse triangulation $\Gamma_H$,
 $$(\tilde I^H w)(y)=\left\{\begin{array}{ll}0, & y\in \partial \Omega,\\
 |\omega_y|^{-1}\int_{\omega_y}w(x)dx, & \hbox{otherwise},
 \end{array}\right.$$
 where $\omega_y$ denotes the union of elements in $\Gamma_H$ that share $y$, $|\omega_y|$ denotes the measure of $\omega_y$.
  Let $I^h$ be the piecewise linear interpolation on the fine triangulation $\Gamma_h$, $Q^h$ be the $L^2$ projector from $\tilde{H}_0^1(\Omega)$ into $V_h$, and
  $\theta_i$ ($i=1,\ldots,J$) be the modified unit partitions of $\{\Omega'_i\}_{i=1}^J$ (see (3.7) in \cite{Toselli1}) which are piecewise linear functions on $\Gamma_h$, satisfying $\sum_{i=1}^J\theta_i=1$, $\max\limits_{1\leq i\leq J}||\nabla \theta_i(x)||_\infty\leq \frac{C}{\delta}$ and $\hbox{supp} (\theta_i)=\Omega'_i$.

  Take $\kappa_0=\tilde I^H$ and define $\kappa_i$, $i=1,2,\ldots,J$ as follows:  for $w\in {\tilde{H}_0^1(\Omega)}$,
$$\kappa_i w:=I^h(\theta_i Q^h(I-\kappa_0)w)+\theta_i(I-Q^h)w.$$
  It is not hard to see that $\kappa_0$ maps $\tilde{H}_0^1(\Omega)$ into $V_0$, $\kappa_i(i=1,\ldots,J)$ maps $\tilde{H}_0^1(\Omega)$ and $V_h$ into $\tilde{H}_0^1(\Omega'_i)$ and $V_i$ respectively, and (\ref{eqn dddd3}) holds. Furthermore, define $\tilde v_i:=I^h(\theta_i Q^h(I-\kappa_0)v), \tilde{\tilde v}_i:=\theta_i(I-Q^h)v,i=1,2,\ldots,J$. It is easy to see that
\begin{equation}\label{eqn dddadda}
\sum\limits_{i=0}^J||\kappa_i v||^2\leq ||\kappa_0 v||^2+C\sum\limits_{i=1}^J||\tilde v_i||^2+C\sum\limits_{i=1}^J||\tilde{\tilde v}_i||^2,
\end{equation}
\begin{equation}\label{eqn dddadd}
\sum\limits_{i=0}^J||\kappa_i v||^2_{\tilde{H}^1(\Omega)}\leq ||\kappa_0 v||^2_{\tilde{H}^1(\Omega)}+C\sum\limits_{i=1}^J||\tilde v_i||^2_{\tilde{H}^1(\Omega)}+C\sum\limits_{i=1}^J||\tilde{\tilde v}_i||^2_{\tilde{H}^1(\Omega)}.
\end{equation}
In order to show that (\ref{eqn ddd10}) and (\ref{eqn ddd9}) hold, it suffices to show that
\begin{equation}\label{eqn dddaddb1}
||\kappa_0 v||^2\leq C||v||^2,
\end{equation}
\begin{equation}\label{eqn dddaddb2}
\sum\limits_{i=1}^J||\tilde v_i||^2\leq C||v||^2,\sum\limits_{i=1}^J||\tilde{\tilde v}_i||^2\leq C||v||^2,
\end{equation}
\begin{equation}\label{eqn dddaddb3}
||\kappa_0 v||^2_{\tilde{H}^1(\Omega)}\leq C||v||^2_{\tilde{H}^1(\Omega)},
\end{equation}
and
\begin{equation}\label{eqn dddaddb4}
\sum\limits_{i=1}^J||\tilde v_i||^2_{\tilde{H}^1(\Omega)}\leq C\left(1+\frac{H}{\delta}\right)||v||^2_{\tilde{H}^1(\Omega)},\quad \sum\limits_{i=1}^J||\tilde{\tilde v}_i||^2_{\tilde{H}^1(\Omega)}\leq C\left(1+\frac{H}{\delta}\right)||v||^2_{\tilde{H}^1(\Omega)}.
\end{equation}
It is known that\cite{Toselli1}
\begin{equation}\label{eqn aad1}
||w-Q^hw||+h||Q^hw||_{\tilde{H}^1(\Omega)}\leq Ch||w||_{\tilde{H}^1(\Omega)}, \quad w\in \tilde{H}_0^1(\Omega).
\end{equation}
By Lemma 3.6 in \cite{Toselli1}, we know
\begin{equation}\label{eqn aad2}
||w-\tilde I^H w||+H||\tilde I^H w||_{\tilde{H}^1(\Omega)}\leq CH||w||_{\tilde{H}^1(\Omega)}, \quad w\in \tilde{H}_0^1(\Omega).
\end{equation}
By a similar proof of Lemma 3.9 in \cite{Toselli1}, we may get
\begin{equation}\label{eqn aad3}
||I^hw||_{\tilde{H}^1(\Omega)}\leq C||w||_{\tilde{H}^1(\Omega)}, \quad   ||I^hw||\leq C||w||,
\end{equation}
here $w$ is a continuous piecewise quadratic function on  $\Gamma_h$.

Then (\ref{eqn dddaddb1}) and (\ref{eqn dddaddb3}) are a direct consequence of (\ref{eqn aad3}), meanwhile (\ref{eqn dddaddb2}) and (\ref{eqn dddaddb4}) can be proved by using (\ref{eqn aad1}), (\ref{eqn aad2}) and (\ref{eqn aad3}) (Please refer to  Lemma 3.12 in \cite{Toselli1} for details).
 \eproof

 Next,  we shall present the main results for the preconditioner system (\ref{preconditioner system}).
 \begin{lemma}\label{lemma n1}
The operator $P_h$ is symmetric positive definite with respect to the inner product $(\cdot,\cdot)_A$.
\end{lemma}
\proof For the proof, we refer to  section 2.3 in \cite{Toselli1}.
\begin{lemma}\label{lemma d4}
For $v\in V_h$, we have
$$(P_hv,v)_A\leq C(v,v)_A.$$
\end{lemma}
\proof
By the definition of $P_0$, we have
\begin{equation}\label{eqn dddd2}
(P_0v,v)_A=(P_0v,P_0v)_A\leq (v,v)_A.
\end{equation} We recall  the index set $I_i$ ($i=1,2,\ldots,J^C$): for any $j,k\in I_i$ with $j\neq k$,
$\Omega'_j\cap \Omega'_k=\emptyset$, so it is easy to see that
\begin{equation}\label{eqn ddd5}
(P_jv,w)_A=(\tilde P_iv,w)_A, \quad \forall w\in V_j, j\in I_i, v\in V_h.
\end{equation}
For any $w\in \tilde V_i$, decompose it as $w=\sum_{j\in I_i}w_j$, $w_j\in V_j$. Then by (\ref{eqn ddd5}), we have
\begin{equation}\label{eqn ddd66}
\sum\limits_{j\in I_i}(P_jv,w_j)_A=\sum\limits_{j\in I_i}(\tilde P_iv,w_j)_A=(\tilde P_iv,w)_A, \quad v\in V_h.
\end{equation}
Taking $w_j=P_jv$ in the above equation, and using Cauchy-Schwarz inequality and  Lemma \ref{lemma d2}, we may obtain
\begin{eqnarray}\label{eqn ddd6}
\sum\limits_{j\in I_i}(P_jv,P_jv)_A&=&(\tilde P_iv,\sum\limits_{j\in I_i}P_jv)_A\nonumber\\
&\leq& ||\tilde P_iv||_A ||\sum\limits_{j\in I_i}P_jv||_A\lesssim ||\tilde P_iv||_A \left(\sum\limits_{j\in, I_i}||P_jv||^2_A\right)^{1/2},
\end{eqnarray}
which is
\begin{equation}\label{eqn ddd7}
\sum\limits_{j\in I_i}(P_jv,P_jv)_A\lesssim ||\tilde P_iv||^2_A.
\end{equation}
It is easy to see that
\begin{equation}\label{eqn ddd8}
(\tilde P_iv,v)_A=(\tilde P_iv,\tilde P_iv)_A\leq (v,v)_A
\end{equation}
and
\begin{equation}\label{eqn dddd1}
\sum\limits_{j\in I_i}(P_jv,P_jv)_A\lesssim ||\tilde P_iv||^2_A\leq (v,v)_A.
\end{equation}
Combining  (\ref{eqn dddd1}) with (\ref{eqn dddd2}), we arrive at
\begin{eqnarray}
\label{eqn ar}(P_hv,v)_A&=&(P_0v,P_0v)_A+(P_1v,P_1v)_A+\cdots+(P_Jv,P_Jv)_A\nonumber\\
&\leq&(P_0v,P_0v)_A+\sum_{j\in I_1}(P_j v,P_jv)_A+\cdots+\sum_{j\in I_{J^C}}(P_jv,P_jv)_A\nonumber\\
&\lesssim& (1+J^C)(v,v)_A.
\end{eqnarray}
Then the lemma is proved.  \eproof
\begin{lemma}\label{lemma d5}
For $v\in V_h$, there exists a decomposition $v=\sum_{i=0}^J v_i,$ $v_i\in V_i$, such that
\begin{equation}\label{eqn dddd11}
\sum\limits_{i=0}^J(v_i,v_i)_A\leq C\left(1+\frac{H}{\delta}\right)||v||^2_A.
\end{equation}
\end{lemma}
\proof Let $\kappa_i$ be the operators as defined in Lemma \ref{lemma d3}, we now prove  the decomposition $v=\sum_{i=0}^J v_i$ with $v_i=\kappa_i v$  may satisfy (\ref{eqn dddd11}).
Since $||w||^2_A\approx ||w||^2_{H^\alpha(\Omega)}$ for $w\in V_h$, it suffices to show that
\begin{equation}\label{eqn dddd22}
\sum\limits_{i=0}^J||\kappa_i v||^2_{H^\alpha(\Omega)}\leq C\left(1+\frac{H}{\delta}\right)||v||^2_{H^\alpha(\Omega)}.
\end{equation}
Actually, using the definition of the norms $||\cdot||_{{\tilde{H}}^{\alpha}(\mathbb{R}^2)}$, $||\cdot||_{\hat{{H}}^{\alpha}(\Omega)}$, Lemma \ref{lemma aaaab} and  Lemma \ref{lemma d3}, we have
\begin{eqnarray}\label{t5}
&&\sum\limits_{i=0}^J||\kappa_i v||^2_{H^\alpha(\Omega)}\approx\sum\limits_{i=0}^J||\kappa_i v||^2_{{\tilde{H}}^\alpha(\mathbb{R}^2)}\nonumber\\
&=&\sum\limits_{i=0}^J\int_0^{\infty}\inf_{w_i\in \tilde{H}^{1}(\mathbb{R}^2)}\left(||\kappa_iv-w_i||^2_{L^2(\mathbb{R}^2)}+t^2||w_i||^2_{\tilde{H}^{1}(\mathbb{R}^2)}\right) t^{-2\alpha-1}dt\nonumber\\
&\leq&\sum\limits_{i=0}^J\int_0^{\infty}\inf_{w_i\in \tilde{H}_0^{1}(\Omega)}\left(||\kappa_iv-w_i||^2_{L^2(\mathbb{R}^2)}+t^2||w_i||^2_{\tilde{H}^{1}(\mathbb{R}^2)}\right) t^{-2\alpha-1}dt\nonumber\\
&=&\sum\limits_{i=0}^J\int_0^{\infty}\inf_{w_i\in \tilde{H}_0^{1}(\Omega)}\left(||\kappa_iv-w_i||^2_{L^2(\Omega)}+t^2||w_i||^2_{\tilde{H}^{1}(\Omega)}\right) t^{-2\alpha-1}dt\nonumber\\
&=&\int_0^{\infty}\inf_{w_0,\ldots,w_J\in \tilde{H}_0^{1}(\Omega)}\sum\limits_{i=0}^J\left(||\kappa_iv-w_i||^2_{L^2(\Omega)}+t^2||w_i||^2_{\tilde{H}^{1}(\Omega)}\right) t^{-2\alpha-1}dt\nonumber\\
&\leq&\int_0^{\infty}\inf_{w\in \tilde{H}_0^{1}(\Omega)}\sum\limits_{i=0}^J\left(||\kappa_iv-\kappa_iw||^2_{L^2(\Omega)}+t^2||\kappa_iw||^2_{\tilde{H}^{1}(\Omega)}\right) t^{-2\alpha-1}dt\nonumber\\
&\lesssim& \left(1+\frac{H}{\delta}\right)\int_0^{\infty}\inf_{w\in \tilde{H}_0^{1}(\Omega)}\left(||v-w||^2_{L^2(\Omega)}+t^2||w||^2_{\tilde{H}^{1}(\Omega)}\right) t^{-2\alpha-1}dt\nonumber\\
&=&\left(1+\frac{H}{\delta}\right)||v||^2_{\hat{{H}}^\alpha(\Omega)}\approx\left(1+\frac{H}{\delta}\right)||v||^2_{\tilde{H}^\alpha(\Omega)}.
  \end{eqnarray}
Moreover,
\begin{eqnarray}\label{t2}||v||^2_{\tilde{H}^{\alpha}(\Omega)}&=&\int_0^{\infty}\inf_{w\in \tilde{H}^{1}(\Omega)}\left(||v-w||^2_{L^{2}(\Omega)}+t^2||w||^2_{\tilde{H}^{1}(\Omega)}\right) t^{-2\alpha-1}dt\nonumber\\
&\leq& \int_0^{\infty}\inf_{w\in \tilde{H}^{1}(\mathbb{R}^2)}\left(||v-w|_\Omega||^2_{L^{2}(\Omega)}+t^2||\,w|_\Omega||^2_{\tilde{H}^{1}(\Omega)}\right) t^{-2\alpha-1}dt\nonumber\\
&\leq& \int_0^{\infty}\inf_{w\in \tilde{H}^{1}(\mathbb{R}^2)}\left(||v-w||^2_{L^{2}(\mathbb{R}^2)}+t^2||w||^2_{\tilde{H}^{1}(\mathbb{R}^2)}\right) t^{-2\alpha-1}dt\nonumber\\
&=&||v||_{\tilde{H}^{\alpha}(\mathbb{R}^2)}\approx ||v||_{H^{\alpha}(\Omega)}.
\end{eqnarray}
Combining (\ref{t2}) with (\ref{t5}), we may get (\ref{eqn dddd22}).\eproof

\begin{lemma}\label{lemma dd1}
For $v\in V_h$, we have
$$(v,v)_A\leq C\left(1+\frac{H}{\delta}\right)(P_hv,v)_A.$$
\end{lemma}
\proof By Lemma \ref{lemma d5}, we know that there exists a decomposition $v=\sum_{i=0}^Jv_i$, $v_i\in V_i$ such that (\ref{eqn dddd11}) holds. Then by the Cauchy-Schwarz inequality and (\ref{eqn dddd11}),  we have
\begin{eqnarray}\label{eqn ee2}
(v,v)_A&=&\sum\limits_{i=0}^J (v,v_i)_A=\sum\limits_{i=0}^J (P_iv,v_i)_A\nonumber\\
&\leq&\left(\sum\limits_{i=0}^J (P_iv,P_iv)_A\right)^{1/2}\left(\sum\limits_{i=0}^J (v_i,v_i)_A\right)^{1/2}\nonumber\\
&=&\left(\sum\limits_{i=0}^J (P_iv,v)_A\right)^{1/2}\left(\sum\limits_{i=0}^J (v_i,v_i)_A\right)^{1/2}\nonumber\\
&\leq&(P_hv,v)^{1/2}_A \left(C\left(1+\frac{H}{\delta}\right)(v,v)_A\right)^{1/2},
\end{eqnarray}
which is Lemma \ref{lemma dd1}\eproof.
\begin{theorem}\label{theorem n} The condition number of $P_h$ satisfies
$$cond(P_h)\leq C\left(1+\frac{H}{\delta}\right).$$
\end{theorem}
\proof The theorem follows from Lemma \ref{lemma d4} and Lemma  \ref{lemma dd1}.\eproof

\section{Numerical results}
\begin{figure}
\begin{center}
  \includegraphics[width=6cm]{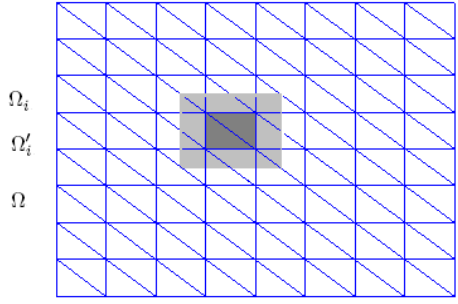}\\
  \caption{The trangulation and the subdomains}\label{fig 3}
\end{center}
\end{figure}

We test PCG method with $B_h$ as a preconditioner.
The tests are carried out by using Matlab software.
The stopping criterion is
 $$||u^k-u^{k-1}||_\infty\leq 10^{-6},$$
where $u^k$ is obtained from $u^{k-1}$ by one step PCG iteration, and  $u^0=0$. Here we only test the case that $\Omega$ is a square domain. We take the triangulations $\Gamma_H,\Gamma_h$ as uniform ones as in figure \ref{fig 3} shows. With viewing the one in figure \ref{fig 3} as the coarse triangulation  $\Gamma_H$,  the subdomains $\Omega_i$ and $\Omega'_i$ are taken as the darker and lighter shaded squares  respectively.

We give two examples in this paper: one is with the probability measure $\tilde{M}$ having a discrete form and the other is with $\tilde{M}$ being a continuous
function. Table \ref{Table1} and Table \ref{Table2}  list the
iterative steps of our PCG methods for Example \ref{exam 1} and Example
\ref{exam 2} respectively. It can be seen that when we fix $\frac{H}{ h}=constant$, for instance $\frac{H}{ h}=8$, the convergence rate of our PCG method
is optimal, i.e., the iterative steps of the PCG method almost keep to be a constant, which coincides with our theoretical results in this paper. From the tables 1-2, we also find that small overlapping case is a good choice for our domain decomposition methods for the SFPDEs, which is same as the overlapping domain decomposition methods for the integer order PDEs.

\begin{example}\label{exam 1} Let $\Omega=[0,2]\times[0,2]$, the equation  is
\begin{equation}\label{problem exam 1}
-\frac14( _{-\infty}D^{1.5}_{x}+  {_{x}}D^{1.5}_{\infty}+
{_{-\infty}D^{1.5}_{y}}+{_{y}D^{1.5}_{\infty}})u=f(x,y),
\end{equation}
where $f(x,y)$ is taken such that the solution is
$u=((2x-x^2)(2y-y^2))^4$,
\end{example}
\begin{table}
  \centering
  \begin{tabular}{||c|c|c|c||}
    \hline
    (h,H)  & $\delta=h$ &  $\delta=2h$  &   $\delta=4h$   \\
    \hline
   $(\frac{2}{64}, \frac{2}{8})$  &  13 & 13 & 14 \\
   \hline
   $(\frac{2}{128}, \frac{2}{16})$  &  12 & 11 & 12 \\
   \hline
   $(\frac{2}{256}, \frac{2}{32})$  &  10 & 10 & 11 \\
   \hline
\end{tabular}
  \caption{The iterative steps of PCG method for Example \ref{exam 1}}\label{Table1}
\end{table}
\begin{table}
  \centering
  \begin{tabular}{||c|c|c|c||}
    \hline
    (h,H)  & $\delta=h$ &  $\delta=2h$  &   $\delta=4h$   \\
    \hline
   $(\frac{2}{64}, \frac{2}{8})$  &  13 & 13& 13 \\
   \hline
   $(\frac{2}{128}, \frac{2}{16})$  &  13 & 13 & 13 \\
   \hline
   $(\frac{2}{256}, \frac{2}{32})$  &  13& 12 & 13 \\
   \hline
\end{tabular}
  \caption{The iterative steps of PCG method for Example \ref{exam 2}}\label{Table2}
\end{table}
\begin{example}\label{exam 2} Let $\Omega=[0,2]\times[0,2]$, the equation is
\begin{equation}\label{problem exam 2}
-D^{1.5}_{\tilde{M}} u=f(x,y),
\end{equation}
where $\tilde{M}(\theta)=1$, $\theta\in [0,2\pi)$, and $f(x,y)$ is taken
as in Example \ref{exam 1}.
\end{example}

\renewcommand{\thesection}{Appendix }
\section{}
\renewcommand{\thesection}{A}
In this appendix, we shall present the definitions of the fractional integrals and derivatives in one
dimension case, and their corresponding properties, which are used in this paper. The related
results can also be found in  \cite{Podlubny1,Samko1}. For function $v(x),x\in
(a,b),$ where $a$ may be taken as $-\infty$, and $b$ as $+\infty$,
the left and right Riemann-Liouville fractional Integrals of order
$\mu>0$ are defined respectively by
\begin{equation}\label{riemann integral}
{_aD_x^{-\mu}}v(x)=\frac1{\Gamma(\mu)}\int_a^x(x-s)^{\mu-1}v(s)ds,\quad
{_xD_b^{-\mu}}v(x)=\frac1{\Gamma(\mu)}\int_x^b(s-x)^{\mu-1}v(s)ds.
\end{equation}
The left and right Riemann-Liouville fractional derivatives of order
$\mu>0$ are defined respectively by
\begin{equation}\label{riemann derivative L}
{_aD_x^{\mu}}v(x)=D^n{_aD_x^{-(n-\mu)}}v(x)=\frac1{\Gamma(n-\mu)}\frac{d^n}{dx^n}\int_a^x(x-s)^{n-\mu-1}v(s)ds,
\end{equation}
\begin{equation}\label{riemann derivative R}
{_xD_b^{\mu}}v(x)=(-D)^n{_xD_b^{-(n-\mu)}}v(x)=\frac{(-1)^n}{\Gamma(n-\mu)}\int_x^b(s-x)^{n-\mu-1}v(s)ds,
\end{equation}
where $n$ is an integer such that $n-1\leq \mu<n$. For $p>0$, $k$ is an integer satisfying $k-1\leq p<k$. If ${_aD_x^p}
v(x)$ is integrable, then
\begin{equation}\label{property 1}
{_aD_x^{-p}}{_aD_x^p}
v(x)=v(x)-\sum\limits_{j=1}^k\left.\left[_aD^{p-j}_x
{v}(x)\right]\right|_{x=a}\frac{(x-a)^{p-j}}{\Gamma(p-j+1)}.
\end{equation}
(\ref{property 1}) is from (2.108) in \cite{Podlubny1}.

\end{document}